\documentclass[]{interact}

\usepackage{epstopdf}
\usepackage[caption=false]{subfig}
\usepackage{url}
\usepackage[numbers,sort&compress]{natbib}
\usepackage{tabularx}
\usepackage{array}
\usepackage{multirow}
\usepackage{hyperref}

\bibpunct[, ]{[}{]}{,}{n}{,}{,}
\renewcommand\bibfont{\fontsize{10}{12}\selectfont}

\theoremstyle{plain}

\theoremstyle{definition}

\theoremstyle{remark}

\begin{document}

\title{A Benchmark of 25 Nonlinear Functions with Domain-Induced Discontinuity-Like for Global Optimization}

\author{
\name{Peicong Cheng\textsuperscript{a}\thanks{CONTACT Peicong Cheng. Email: cheng.p.b2bf@m.isct.ac.jp} and Makoto Yamashita\textsuperscript{a}}
\affil{\textsuperscript{a}Department of Mathematical and Computing Science, Institute of Science Tokyo, Tokyo, Japan}
}

\maketitle

\begin{abstract}
A benchmark of 25 nonlinear optimization problems exhibiting domain-induced discontinuity-like behavior is proposed to support the performance evaluation of global optimization algorithms under feasibility-scarce and structurally disconnected landscapes. Referred to as the CPC Benchmark (\textit{Challenging Problems for Computation}), the test suite consists of functions that are continuous on their natural domains, while infeasible regions and undefined evaluations are implicitly embedded in the objective, creating substantial challenges for global minimization. Six representative algorithms from diverse methodological paradigms are assessed to investigate the structural difficulty and discriminative capability of the proposed benchmark. Numerical results show that many functions possess extremely small feasible regions and strong precision sensitivity near feasibility boundaries, complicating initialization, feasibility discovery, and reliable objective evaluation. The findings demonstrate that the CPC benchmark provides a rigorous and discriminative testbed for advancing research on global optimization under hidden feasibility constraints.\end{abstract}


\section{Introduction}\label{sec:section-1-Introduction}

Benchmark functions are fundamental tools for the performance evaluation of global optimization algorithms and play a central role in the interface between algorithm design and software implementation. While theoretical analysis establishes algorithmic soundness, empirical evaluation provides critical insights into practical performance, including robustness, global optimality, and behavior under diverse conditions. Moreover, the development of challenging benchmark problems not only complements theoretical advances but also drives the evolution of optimization algorithms. In practice, meaningful progress in algorithm design is often achieved by tackling difficult optimization tasks.

However, despite the rapid development of new global optimization algorithms in recent years, many studies still rely on benchmark functions introduced more than a decade ago for performance evaluation. For instance, recent high-impact algorithms such as Harris Hawks Optimization \cite{Heidari2019} (2019) and Marine Predators Algorithm \cite{Faramarzi2020} (2020) have attracted considerable attention and demonstrated promising performance, yet their numerical evaluations are primarily based on classical benchmarks \cite{Yao1999, Digalakis2001} developed in the early 2000s. This reflects a broader trend in which outdated benchmark suites are repeatedly used, limiting not only the reliability and discriminative power of algorithmic comparisons but also the development of modern software infrastructures for global optimization research. As a result, many of these benchmark functions have become insufficiently challenging, as they can be efficiently solved by modern optimization algorithms or professional mathematical solvers. Therefore, the development of more challenging and representative benchmark problems is essential for a fair and comprehensive assessment of modern optimization algorithms.

Furthermore, although discontinuous objective functions are common in various real-world applications, very few benchmark suites have been designed to explicitly capture this characteristic. Existing studies on discontinuous optimization are therefore often constrained by the lack of appropriate test problems. For example, Ban and Yamazaki \cite{Ban2019} proposed a global optimization method tailored for discontinuous problems with infeasible regions. However, due to the lack of standardized benchmark functions of this type, their evaluation relied on continuous benchmark functions with artificially constructed  infeasible region. This workaround highlights a fundamental gap: current benchmark collections do not adequately represent optimization problems exhibiting domain-induced discontinuity-like behavior , where infeasible regions are inherently embedded in the function domain rather than externally imposed.

In this study, we focus on a class of minimization problems exhibiting domain-induced discontinuity-like behavior (DIDL behavior). Such problems can be interpreted as implicitly constrained optimization problems, where the objective function remains continuous on its natural domain but contains undefined or infeasible regions within the prescribed search space. To address the lack of standardized benchmarks for this class of problems, we propose a new benchmark suite consisting of 25 challenging nonlinear optimization problems with DIDL behavior, referred to as the CPC Benchmark (Challenging Problems for Computation). The CPC benchmark is designed not only as a collection of challenging test problems but also as a standardized evaluation framework for assessing the robustness and effectiveness of global optimization algorithms. The proposed functions exhibit diverse structural properties and substantially greater difficulty than many existing problems exhibiting discontinuity-related behavior.

The rest of this paper is organized as follows: Section~\ref{sec:section-2-methodology} presents the methodology, including problem definition, taxonomy of discontinuities, algorithm selection, and evaluation protocol. Section~\ref{sec:section-3-characteristics} analyzes six key characteristics of optimization problems with DIDL behavior and summarizes the structural properties of all proposed benchmark functions. Section~\ref{sec:section-4-numerical-experiment} reports experimental setup and numerical results. Section~\ref{sec:section-5-Discussion} provides a structured discussion of the benchmark from both theoretical and algorithmic perspectives. Finally, Section~\ref{sec:section-6-conclusion} concludes the paper with remarks and future research directions.

\section{Methodology}\label{sec:section-2-methodology}

\subsection{Problem Definition}\label{subsec:problem-definition}

We consider the following global minimization problem:

\[
\min f(\mathbf{x})
\quad \text{s.t.} \quad
\mathbf{x}\in\mathcal I,
\]

where $\mathcal I\subset\mathbb R^n$ denotes the prescribed search domain. The objective function is defined as

\[
f:\mathcal D \rightarrow \mathbb R,
\]

where $\mathcal D\subseteq\mathcal I$ represents the natural domain of the function. In general, $\mathcal D$ is not explicitly provided to the optimization algorithm and may constitute only a fragmented subset of the prescribed search domain $\mathcal I$.

In this study, the observed optimization difficulty arises not from explicitly imposed infeasible regions or artificially introduced discontinuities, but from intrinsic domain restrictions embedded in the analytical form of the objective function. Elementary functions such as logarithmic ($\log(x)$) and square roots ($\sqrt{x}$) possess intrinsic domain restrictions. Consequently, when the prescribed search domain is not fully contained within the natural domain of the objective function, infeasible or undefined regions emerge within the search space. As a result, the effective feasible region may become fragmented, producing discontinuity-like behavior from the viewpoint of optimization algorithms, even though the objective function remains continuous on its natural domain.

Throughout this study, unless otherwise stated, the admissible domain refers to the computationally admissible domain, namely the set of points at which the objective function can be successfully evaluated using standard real-valued floating-point arithmetic. For example, the expression $(-8)^{1/3}$ is mathematically equal to $-2$; However, when evaluated using standard floating-point power operators in numerical computing environments, the same expression may return $\mathrm{NaN}$ or trigger a domain error depending on the programming implementation \cite{ISO9899_2024}. Therefore, the admissible domain adopted in this study is defined in terms of computational evaluability rather than purely mathematical admissibility.

To better characterize the origins of DIDL behavior in the proposed benchmark functions, we introduce and analyze five representative nonlinear operator whose domain restrictions may induce fragmented feasible regions within the prescribed search space. The domains listed below correspond to computationally admissible domains under standard real-valued floating-point arithmetic.

\begin{itemize}

    \item \textbf{Logarithmic Function:} $f(x) = \log(x)$, where $\text{dom}(f) = \mathbb{R}_{++}$.
    
    \item \textbf{Reciprocal Function:} $f(x) = \frac{1}{x}$, where $\text{dom}(f) = \mathbb{R} \setminus \{0\}$.
    
    \item \textbf{Square Root Function:} $f(x) = \sqrt{x}$, where $\text{dom}(f) = \mathbb{R}_{+}$.
    
    \item \textbf{General Power Function:} $f(x) = x^p$, where $\text{dom}(f)$ depends on $p$:
    \begin{itemize}
    \item If $p \in \mathbb{Z}_{\geq 0}$, then $\text{dom}(f) = \mathbb{R}$.
    \item If $p \in \mathbb{Z}_{< 0}$, then $\text{dom}(f) = \mathbb{R} \setminus \{0\}$.
    \item If $p \in p \in (\mathbb{R} \setminus \mathbb{Z})_{\ge 0}$, then $\text{dom}(f) = \mathbb{R}_{+}$.
    \item If $p \in (\mathbb{R} \setminus \mathbb{Z})_{< 0}$, then $\text{dom}(f) = \mathbb{R}_{++}$.
    \end{itemize}
    
    \item \textbf{Exponential Function:} $f(x) = a^x$, where $\text{dom}(f)$ depends on $a$:
    \begin{itemize}
    \item If $a \in \mathbb{R}_{++}$, then $\text{dom}(f) = \mathbb{R}$.
    \item If $a := 0$, then $\text{dom}(f) = \mathbb{R}_{+}$.
    \item If $a \in \mathbb{R}_{--}$, then $\text{dom}(f) = \mathbb{Z}$.
    \end{itemize}
\end{itemize}

The reciprocal function $\frac{1}{x}$ and the square root function $\sqrt{x}$ are special cases of the general power function, corresponding to $p=-1$ and $p=\tfrac12$, respectively. Owing to their frequent use in practical modeling and optimization problems, they are discussed separately.

\begin{figure}
    \centering
    \includegraphics[width=0.6\linewidth]{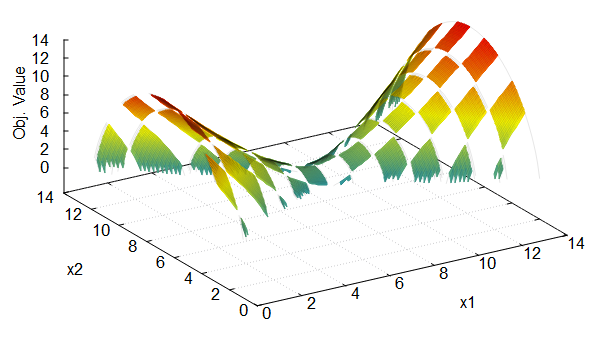}
    \caption{Function graph of CPC-DF8.}
    \label{fig:cpc-demo-df8}
\end{figure}

\begin{figure}
    \centering
    \includegraphics[width=0.6\linewidth]{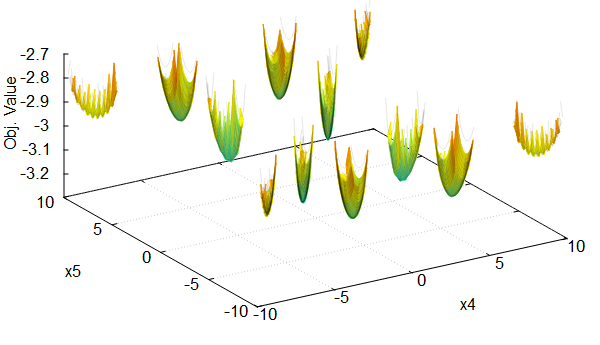}
    \caption{Function graph of CPC-DF12.}
    \label{fig:cpc-demo-df12}
\end{figure}

Objective functions composed of combinations of these elementary nonlinear functions may exhibit a highly complex optimization landscape, potentially resulting in multiple disconnected discontinuity-like feasible regions. More importantly, when these nonlinear operators involve multiple decision variables, the resulting computationally admissible domain is often determined jointly by several variables rather than independently for each variable, giving rise to complex hidden feasibility structures. To provide a clearer illustration of DIDL behavior, two representative functions from the CPC Benchmark, CPC-DF8 and CPC-DF12, are presented below together with their corresponding function graphs in Fig.~\ref{fig:cpc-demo-df8} and Fig.~\ref{fig:cpc-demo-df12}.

Representative Example 1: CPC-DF8

\begin{equation*}   
    \begin{aligned}
        f(\textbf{x}) := &{1- \left( \frac{ \mathrm{sin} \left( \pi \cdot \left( x_{1}-2 \right) \right) \cdot \mathrm{sin} \left( \pi \cdot \left( x_{2}-2 \right) \right) }{\pi^{2} \cdot x_{1} \cdot \left( x_{1}-2 \right) \cdot \left( x_{2}-2 \right) } \right) ^{1.03} + } \\
        & {\left( 2+ \left( x_{1}-7 \right) ^{2}-2 \cdot \left( x_{2}-7 \right) ^{2} \right) ^{0.65}}
    \end{aligned}
    \tag{CPC-DF8}
    \label{eq:cpc-df8}
\end{equation*}

where $x_i \in [0, 14], \quad i = 1, 2$

Representative Example 2: CPC-DF12
\begin{equation*}   
    \begin{aligned}
        &f(\textbf{x}) := {-\sum\limits _{i=1}^{n-1=4} \left( \frac{ \sqrt{ \sin \left( \sqrt{x_{i+1}^{2}-x_{i}^{2}}-0.5 \right) -0.5}}{ \left( 10 \cdot \left( x_{i+1}^{2}+x_{i}^{2} \right) -0.85 \right) ^{0.2}}+0.5 \right) }
    \end{aligned}
    \tag{CPC-DF12}
\end{equation*}

Where $x_i \in [-100, 100], \quad i = 1, \ldots, n$

As illustrated by the two representative examples above, the resulting function landscapes exhibit highly fragmented feasible regions and disconnected search-space structures within the prescribed domain. These characteristics arise from intrinsic domain restrictions embedded in the objective function, which implicitly exclude portions of the search space by rendering them infeasible or undefined.

\subsection{Taxonomy of Discontinuity-Related Behaviors}

Discontinuity in global optimization is a broad and often loosely defined concept. In existing benchmark collections, such as those surveyed in \cite{Jamil2013Survey}, many problems are labeled as "discontinuous," although the underlying sources of discontinuity differ significantly. To reduce this ambiguity, we classify discontinuity-related behaviors into four representative categories based on their structural origin.

It is important to emphasize that not all these categories correspond to strictly defined mathematical discontinuities. Some functions remain continuous over their natural domain but exhibit numerical or structural characteristics that cause optimization algorithms to behave as if discontinuities were present. In other cases, functions are continuous on their natural domain but appear effectively discontinuous over the prescribed search domain due to regions where the analytical expression is undefined.

\begin{figure}
    \centering
    \includegraphics[width=0.445\textwidth]{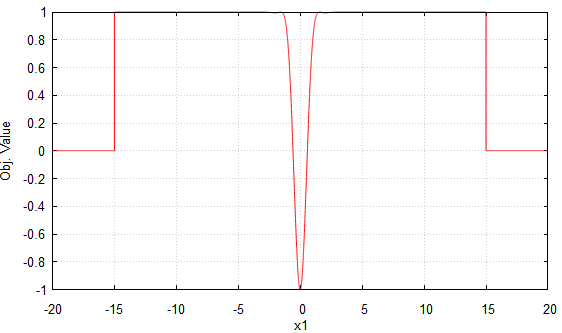}
    \includegraphics[width=0.445\textwidth]{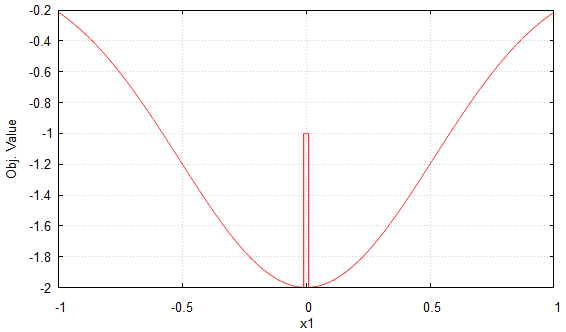}
    \caption{Illustration of numerically induced pseudo-discontinuity in the Xin-She Yang Function 3 (F170). The global view on the interval [-20,20] (left) shows an apparently flat plateau with abrupt value transitions around $x = \pm15$, whereas the local zoom on [-1,1] (right) reveals a sharply curved and narrowly concentrated basin around the origin. These sharply localized variations cause the function to appear discontinuous under finite-precision sampling, despite being mathematically continuous over its domain.}
    \label{fig:f170-pseudo-discontinuity}
\end{figure}

\begin{figure}
    \centering
    \includegraphics[width=0.6\linewidth]{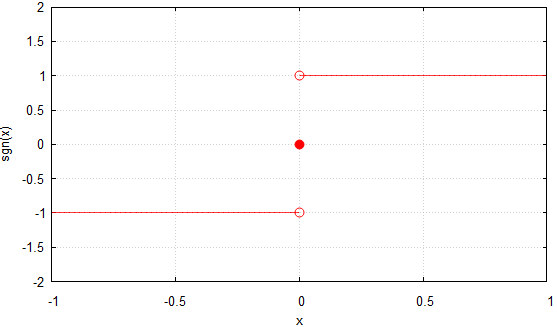}
    \caption{The sign function $f(x)=\mathrm{sgn}(x)$, showing a jump discontinuity at the origin arising from its piecewise-defined structure.}
    \label{fig:sign-piecewise-discontinuity}
\end{figure}

\begin{figure}
    \centering
    \includegraphics[width=0.6\linewidth]{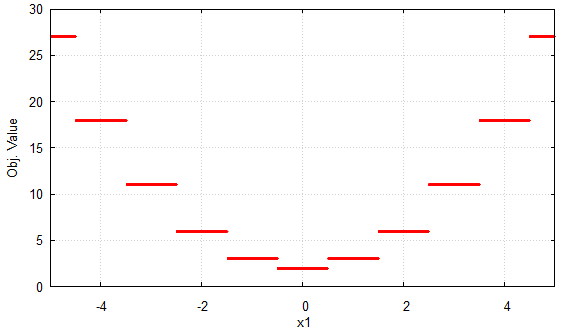}
    \caption{Illustration of the Step Function 2. The objective values remain constant within each quantized interval of $x_1$, forming flat plateaus, and exhibit abrupt jumps when transitioning between adjacent levels.}
    \label{fig:stepwise-quantization}
\end{figure}

\begin{enumerate}
    \item \textbf{Numerically Induced Pseudo-Discontinuity}: This category refers to functions that are mathematically continuous over the entire domain but exhibit extremely steep gradients, sharp curvature variations, or rapid oscillations within very small neighborhoods. As a result, small changes in the input may lead to large variations in function values at the numerical scale, causing optimization algorithms to perceive discontinuous behavior.

    A representative example is the Xin-She Yang Function 3 (F170 in \cite{Jamil2013Survey}). Although the function is continuous over its domain, its landscape exhibits sharp localized variations and extreme value changes within very narrow regions. As illustrated in Fig.~\ref{fig:f170-pseudo-discontinuity}, the function exhibits steep gradients and abrupt value changes at several locations, which may cause it to appear numerically discontinuous.
    
    These steep gradients and localized structures do not introduce true mathematical discontinuities, but they create significant numerical sensitivity. Consequently, small perturbations in the input can lead to large variations in function values, causing optimization algorithms to exhibit unstable search behavior and perceive the landscape as effectively fragmented. 

    \item \textbf{Piecewise Structural Discontinuity}: This category comprises functions that are explicitly defined by distinct analytical expressions over multiple subdomains. Discontinuities arise at the interfaces between these pieces whenever the function values or derivatives fail to match. Such functions are intrinsically discontinuous on their natural domain in the strict mathematical sense. Representative examples include the sign function $f(x)=\mathrm{sgn}(x)$, which exhibits a jump discontinuity at the origin. An illustration is provided in Fig.~\ref{fig:sign-piecewise-discontinuity}.

    \item \textbf{Stepwise Quantization Discontinuity}: This category characterizes functions whose outputs are quantized into discrete stepwise levels, producing flat plateaus separated by abrupt jumps. Such functions are also intrinsically discontinuous on their natural domain in the strict mathematical sense. Typical examples include the floor operator $\lfloor x \rfloor$ and the ceiling operator $\lceil x \rceil$. A representative benchmark example is the Step Function 2 (F138 in \cite{Jamil2013Survey}). An illustration is provided in Fig.~\ref{fig:stepwise-quantization}.

    \item \textbf{Domain-induced Discontinuity-like Behavior (DIDL Behavior)}: This category refers to optimization problems whose objective functions remain continuous on their natural domains but possess intrinsic domain restrictions arising from mathematical operators such as fractional powers, logarithmic terms, or exponential functions introduced in Section~\ref{subsec:problem-definition}. When these restrictions are violated, the objective function becomes undefined and cannot be evaluated. Under a prescribed search domain, this resulting exclusion of infeasible regions may fragment the admissible search space into multiple disconnected components. Consequently, optimization algorithms may encounter abrupt transitions between feasible and infeasible regions, causing the search landscape to exhibit discontinuity-like behavior. From a mathematical perspective, the objective function remains continuous when restricted to its natural domain. Therefore, the observed behavior does not originate from a true discontinuity of the function itself, but from the fragmentation of the admissible domain induced by embedded feasibility conditions. All functions in the proposed CPC Benchmark belong to this category, as illustrated in Fig.~\ref{fig:cpc-demo-df8} and Fig.~\ref{fig:cpc-demo-df12}.
\end{enumerate}

Among these four categories, DIDL behavior is the only category that fundamentally alters the feasible search space. For the other three categories, all points within the prescribed domain remain mathematically well-defined and evaluable. 

To highlight these distinctions from both mathematical and optimization perspectives, Table~\ref{tab:discontinuity-comparison} summarizes the key properties of the four discontinuity types.

\begin{table}
\tbl{Comparative characteristics of four discontinuity-related categories}
{\begin{tabularx}{\textwidth}{
    >{\centering\arraybackslash}p{4.5cm}
    >{\centering\arraybackslash}X
    >{\centering\arraybackslash}X}
\toprule
Discontinuity Type & Mathematical Continuity & Feasibility Restriction\\
\midrule
Numerically induced pseudo discontinuity & Continuous & 
No \\ 
Piecewise discontinuity & Typically discontinuous & No\\
Stepwise discontinuity & Discontinuous & No\\
Domain-Induced discontinuity-like behavior & Continuous & Yes\\
\bottomrule
\end{tabularx}}
\label{tab:discontinuity-comparison}
\end{table}

\subsection{Overview of Selected Algorithms}\label{subsec:overview-of-selected-algorithms}

To comprehensively evaluate the CPC Benchmark, six representative global optimization algorithms (GOAs) were selected to cover three major methodological dimensions: (i) deterministic and stochastic; (ii) gradient-based and gradient-free; (iii) heuristic and exact. These algorithms are well-established and widely recognized in the global optimization literature. 

A consolidated overview and classification of the six algorithms is provided in Table~\ref{tab:algorithm-overview}. Detailed descriptions of their computational platforms and parameter configurations are presented below. Because algorithmic performance may vary across computational platforms, we adopt the most established and representative platform for each algorithm to ensure a fair comparison.

\subsubsection*{Branch \& Bound (B\&B)}

\textbf{Introduction:} Branch \& Bound (B\&B)\cite{Horst2013Book} is a classic deterministic, gradient-free, exact algorithm designed primarily for discrete and combinatorial optimization, especially integer programming problem. It systematically explores subsets of the solution space using upper and lower bounds to prune unpromising branches. Moreover, B\&B has also demonstrated strong capabilities in nonlinear global optimization \cite{Misener2014, Amaran2012, Wanufelle2007}. B\&B is theoretically guaranteed to converge to a global optimum under suitable assumptions, including finite-dimensionality, bounded feasible regions, continuity of the function, and valid bounding schemes that become tight as partitions are refined~\cite{Floudas2013Book, Kannan2018}. Owing to this convergence guarantee, B\&B has been widely adopted in professional optimization application.

\textbf{Algorithm Platform:} LINGO~\cite{LINGO2025} is a global optimization solver whose global solver toolbox is primarily based on the B\&B framework. LINGO is a widely recognized and reliable commercial optimization software that demonstrates strong performance in solving nonlinear global optimization problems. Although BARON~\cite{BARON2026} is widely recognized as one of the leading B\&B-based global optimization solver, it was not included in this study due to its limitation in handling trigonometric functions such as $\sin$ and $\cos$. Since most benchmark functions proposed in this study contain trigonometric terms, BARON cannot be applied directly. Therefore, LINGO was selected as the representative B\&B solver. LINGO 20.0 (academic license) was used in this study. 

\textbf{Parameter Setting:} The precision parameter was set to 16. Other parameters remained at default.

\subsubsection*{Universal Global Optimization (UGO)}

\textbf{Introduction:} UGO is a proprietary stochastic global optimization algorithm developed within the 1stOpt software~\cite{FirstOpt2024}. Although its internal mechanism has not been disclosed, the algorithm consistently demonstrates strong performance on small-scale nonlinear global optimization problems. A distinctive feature of UGO is that it does not necessarily require user-specified initial points or variable bounds; even without such information, it is often able to rapidly identify high-quality, and sometimes unexpectedly good solutions. This robustness highlights the intrinsic strength of UGO's search mechanism.

\textbf{Algorithm Platform:} 1stOpt 11.0.

\textbf{Parameter Setting:} The default settings of UGO in 1stOpt were used.

\subsubsection*{Genetic Algorithm (GA)}

\textbf{Introduction:} GA is a classic population-based stochastic metaheuristic algorithm, and is inspired by the principles of natural selection and genetics. It operates through evolutionary processes including selection, crossover, and mutation to explore the solution space\cite{Holland1975}\cite{Goldberg1989}. The gradient-free nature of GA makes it particularly well-suited for solving non-convex global optimization problems. Due to its strong global search capability and robustness against local minima, GA has been widely recognized as a standard and effective method in global optimization\cite{Norouzi2014}.

\textbf{Algorithm Platform:} MATLAB 2024b. The built-in \texttt{ga} function in MATLAB was used for implementing GA.

\textbf{Parameter Setting:} Population size and maximum generations were both set to 2000. Initial population shuffling was enabled to ensure stochastic initialization. Other parameters remained at MATLAB's default values.

\subsubsection*{Simulated Annealing (SA)}

\textbf{Introduction:} SA is a typical stochastic, single-solution-based metaheuristic algorithm inspired by the annealing process in metallurgy. It accepts worse solutions at high temperatures to enhance global exploration and gradually reduces this acceptance as the temperature decreases\cite{Kirkpatrick1983}, enabling convergence toward high-quality solutions while avoiding local optima\cite{Locatelli2000}.Due to its simplicity, flexibility, independent from gradient information, and strong global search capability, SA has become a widely used and effective method for global optimization\cite{Aarts1989}.

\textbf{Algorithm Platform:} MATLAB 2024b. The built-in \texttt{simulannealbnd} function in MATLAB was used for implementing the SA.

\textbf{Parameter Setting:} Maximum iteration were set to 2000. Initial population shuffling was enabled to ensure stochastic initialization. Other parameters remained at default values.

\subsubsection*{A quasi-Newton method proposed by Broyden, Fletcher, Goldfarb, and Shanno (BFGS) + Mutil-Start}

\textbf{Introduction:} BFGS is a deterministic, gradient-based local optimization method from the quasi-Newton family. By iteratively approximating the inverse Hessian matrix, BFGS achieves fast convergence near local minima without requiring explicit second-order derivatives\cite{Broyden1970, Fletcher1970, Goldfarb1970, Shanno1970}. Although inherently local, combining BFGS with a multi-start strategy allows the algorithm to explore multiple regions of the search space\cite{Tsoulos2020,Nocedal1999Book}, thereby serving as an effective heuristic approach for global optimization \cite{Puthuparampil2017}.

\textbf{Algorithm Platform:} MATLAB 2024b. The built-in \texttt{fminunc} function with the \texttt{quasi-newton} algorithm was used, and MATLAB's \texttt{MultiStart} solver was employed to perform multi-start optimization.

\textbf{Parameter Setting:} The maximum iterations were set to 2000. In each run, 2000 randomized initial points were generated within the given search range and used to start the BFGS algorithm iteration. Other parameters remain at default values.

\subsubsection*{Gradient-based Particle Swarm Optimization Algorithm (GPSO)}

\textbf{Introduction:} GPSO is a hybrid metaheuristic that integrates the global exploration ability of Particle Swarm Optimization (PSO) with local refinement through gradient-based updates\cite{Noel2012}. By augmenting particle movements with gradient information, GPSO accelerates convergence and improves solution precision while maintaining PSO's global search behavior, making it an effective and competitive approach for solving global optimization problems\cite{Zulu2023}.

\textbf{Algorithm Platform:} Matlab 2024b + PlatEMO \cite{Tian2017PlatEMO}. Since there's no built-in function for GPSO in Matlab, PlatEMO platform was used to implement GPSO through the command \texttt{platemo('algorithm', @GPSO)}.

\textbf{Parameter Setting:} The population size was set to 2000, and the maximum number of function evaluations (\texttt{maxEF}) was set to 100000. Other parameters remain at default values.\par

\begin{table}
    \tbl{Overview of six representative global optimization algorithms}
    {
        \begin{tabularx}{\textwidth}{
            >{\centering\arraybackslash}p{2.3cm}
            >{\centering\arraybackslash}p{2.5cm}
            >{\centering\arraybackslash}p{1.8cm}
            >{\centering\arraybackslash}X
            >{\centering\arraybackslash}p{2.8cm}}
        \toprule
        Global Algorithm & Platform & Deterministic/ Stochastic & Gradient-based/ Gradient-free & Heuristic/ Exact \\
        \midrule
        B\&B & LINGO 20.0 & Deterministic & Gradient-free & Exact \\
        UGO & 1stOpt 11.0 & Stochastic & - & - \\
        GA & MATLAB 2024b& Stochastic & Gradient-free & Heuristic (multi-population) \\
        SA & MATLAB 2024b& Stochastic & Gradient-free & Heuristic (single-solution) \\
        BFGS+MS & MATLAB 2024b& Deterministic & Gradient-based & Heuristic \\
        GPSO & MATLAB 2024b + PlatEMO & Stochastic & Gradient-based & Heuristic \\
        \bottomrule
    \end{tabularx}
    }
    \tabnote{B\&B: Branch \& Bound; UGO: Universal Global Optimization; GA: Genetic Algorithm; SA: Simulated Annealing; BFGS: Broyden–Fletcher–Goldfarb–Shanno method; MS: Multi-Start; GPSO: Gradient-based Particle Swarm Optimization.}
    \label{tab:algorithm-overview}
\end{table}

\subsection{Evaluation Protocol}\label{subsec:evaluation-protocol}

We conducted all the numerical experiments on Windows 11 with Intel Core i7-12700H and 32 GB memory space. All computational experiments were conducted under the same hardware environment to ensure fairness, and all numerical computations were performed in double precision.

In the numerical experiments, each algorithm was executed 20 times on every benchmark function, and the best result among these runs was recorded.

Because MATLAB will return complex values when evaluating points in undefined or infeasible regions, a penalty function was applied only to the MATLAB-based algorithms. Since our objective is to obtain real-valued solutions, any point falling in an infeasible region was assigned a penalty value of $10^{100}$, thereby removing such regions from the feasible search space. To illustrate this procedure, we provide a concrete example based on the CPC-DF13.

The original mathematical definition of CPC-DF13 is given as follows:
\begin{equation}
    \begin{aligned}
        f(\textbf{x}) := &{ \sqrt{x_{1}-2 \cdot x_{2}^{2}- \exp \left( x_{2}-x_{1}^{2} \right) }-}\\
        &{\left( 0.5 \cdot \cos \left( 3 \cdot \pi \cdot x_{1}+4 \cdot \pi \cdot x_{2}+5 \right) -0.495 \cdot x_{1} \right) ^{\frac{1}{5}}}
    \end{aligned}
    \tag{CPC-DF13}
    \label{eq:cpc-df13}
\end{equation}

Where $x_i \in [-100, 100], \quad i = 1, 2$

To ensure that MATLAB returns only real-valued outputs, the objective function is modified with a penalty mechanism, and the penalized version of CPC-DF13 is given in the piecewise definition below.

\begin{equation*}
    f(\mathbf{x}) :=
        \begin{cases}
            10^{100}, & \text{if } base_{13\text{-}1} < 0, \\
            10^{100}, & \text{if } base_{13\text{-}2} < 0, \\
            \sqrt{base_{13\text{-}1}} - base_{13\text{-}2}^{\frac{1}{5}}, & \text{otherwise}.
        \end{cases}
\end{equation*}

Where: \begin{equation*}
    \begin{aligned}
        base_{13\text{-}1} &:= x_{1}-2 \cdot x_{2}^{2}- \exp \left( x_{2}-x_{1}^{2} \right), \\
        base_{13\text{-}2} &:= 0.5 \cdot \cos \left( 3 \cdot \pi \cdot x_{1}+4 \cdot \pi \cdot x_{2}+5 \right) -0.495 \cdot x_{1}.
    \end{aligned}
\end{equation*}

\section{Structural Characteristics of Domain-Induced Discontinuity-Like Problems}\label{sec:section-3-characteristics}

Optimization problems with DIDL behavior possess several distinctive structural characteristics that fundamentally influence the search process and optimization performance. To systematically characterize these properties in the proposed CPC benchmark, we analyze the benchmark functions from six complementary perspectives: feasible region structure, modality, multiplicity of optimal solutions, separability, dimensionality, and solution sensitivity.

\subsection{Feasible Region Structure}

Feasible region of a function describes the subset of the search space in which the function is mathematically well-defined and returns real-valued outputs. Correspondingly, the feasible ratio is defined as:
\begin{align*}
\textrm{Feasible Ratio (theoretical)} := \frac{\textrm{Volume of the feasible region}}{\textrm{Volume of the entire search space}}
\end{align*}

For optimization problems whose prescribed search domain is entirely contained within the natural domain of the objective function, the feasible ratio is 100\%, since every point in the search space is evaluable. In contrast, optimization problems with DIDL behavior typically possess substantially smaller feasible ratios due to intrinsic domain restrictions embedded in the objective function.

A low feasible ratio substantially increases the difficulty of locating valid search regions, as large portions of the search domain may correspond to undefined or infeasible points. Consequently, global optimization algorithms may struggle to generate admissible initial solutions, and in extreme cases, may fail to identify any feasible point, leading to premature termination or invalid outputs.

In practice, the exact measure of the feasible region is generally intractable for problems with DIDL behavior. Therefore, we employ the Monte Carlo method to approximate the feasible ratio:
\begin{align*}
\textrm{Feasible Ratio (Monte Carlo approximation)} := \frac{\textrm{Number of sampled feasible points}}{\textrm{Total number of sampled points}}
\end{align*}

In particular, for each function, we draw five million random samples within the search space per run, repeat this procedure five times, and take the average feasible ratio across the five runs as the final estimate.

\begin{figure}
    \centering
    \includegraphics[width=0.445\textwidth]{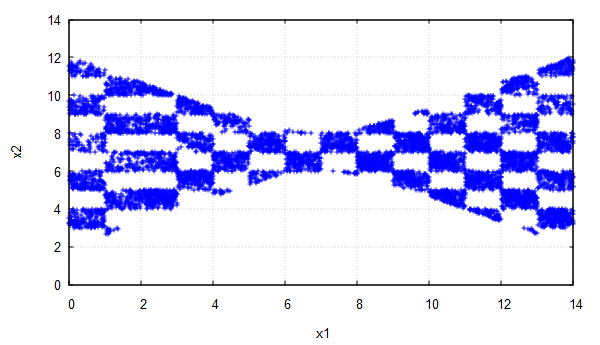}
    \includegraphics[width=0.445\textwidth]{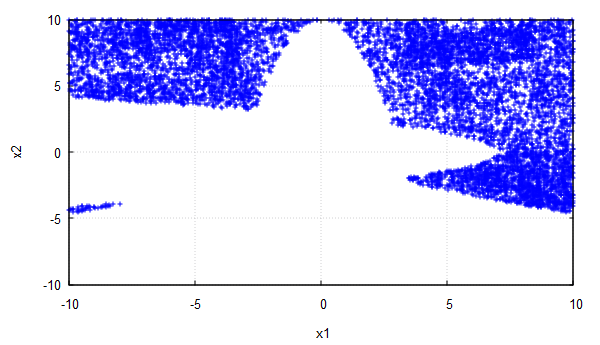}
    \caption{Illustration of the feasible region for two benchmark functions: CPC-DF8 (left) and CPC-DF9 (right). Feasible samples are shown in blue.}
    \label{fig:feasible_regions}
\end{figure}

Figure~\ref{fig:feasible_regions} provides visual examples of the feasible regions for two representative benchmark functions, illustrating the spatial sparsity and fragmentation commonly observed in problem with DIDL behavior.

\subsection{Dimensionality}

Dimensionality refers to the number of decision variables in an optimization problem. Some benchmark functions are defined with a fixed dimensionality (e.g., 2 or 4), while others are formulated in a scalable manner, where the dimension is denoted symbolically (e.g. $n=3$, $n=4$). Such scalable functions can be extended to arbitrary dimensions, allowing researchers to evaluate algorithmic performance under varying problem size.

In general, optimization problems become more challenging as dimensionality increase, particularly when the objective function is non-separable. To support research on large-scale optimization, many benchmark suites\cite{Xu2023LS,Qiao2024LSCMO} provide scalable test functions. These functions are typically designed to assess an algorithm's behavior as the dimension grows. However, their low-dimensional instances often exhibit relative simple structures - frequently fully separable - and can be solved to global optimality with ease by professional mathematical solvers. As a result, they primarily evaluate scalability rather than an algorithm's ability to precisely identify global optima in challenging landscapes.

Some large-scale benchmarks~\cite{Plevris2022} also place the global optimum at a trivial location such as $(0,0,\dots,0)$, which can result in misleading performance evaluations, since algorithms that initialize at or near the origin may reach the optimum without demonstrating effective search capability.

In contrast, the CPC benchmark focuses primarily on low-dimensional yet structurally challenging problems. Rather than competing with existing large-scale benchmark suites, the CPC benchmark provides a complementary perspective by emphasizing difficulty arising from domain-induced discontinuity-like behavior, fragmented feasible regions, and hidden feasibility restrictions. Consequently, it offers a test environment that highlights an algorithm's robustness in locating high-quality feasible solutions within structurally complex search spaces.

\subsection{Separability}

Separability describes whether a multivariate objective function can be expressed as the sum of independent single-variable components. Formally, a function $f(\mathbf{x})$ is said to be \emph{separable} if it can be written as 
\[
f(\mathbf{x}) = \sum_{i=1}^{n} f_i(x_i).
\]

In separable functions, the optimization problem can be reduced to a set of one-dimensional problems that can be solved independently. This property significantly simplifies the optimization process and often allows algorithms to perform more efficiently.

In contrast, non-separable functions contain cross terms among variables, meaning that changes in one variable affect the behavior of the entire function. This increases the complexity and poses a greater challenge for optimization algorithms.

In optimization problems with DIDL behavior, the role of separability becomes even more important because feasibility is governed by implicit domain restriction. For separable functions, feasibility can often be assessed independently for each variable, reducing the identification of infeasible regions to a collection of one-dimensional conditions. In contrast, for non-separable functions, feasibility depends on the joint configuration of multiple variables. As a result, infeasible regions arise from interactions among variables and often form highly irregular and fragmented structures. This substantially increase the difficulty of feasibility detection and search-space exploration, as optimization algorithms must navigate complex high-dimensional boundaries between feasible and infeasible regions.

\begin{figure}
    \centering
    \includegraphics[width=0.45\textwidth]{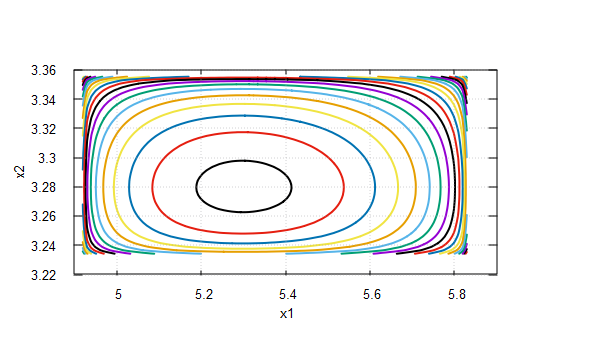}
    \includegraphics[width=0.45\textwidth]{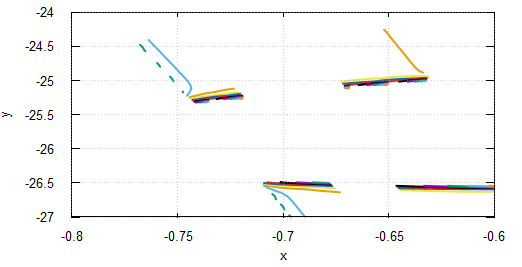}
    \caption{Visualization of separability structures for CPC-DF7 (left) and CPC-DF22 (right). The separable function CPC-DF7 exhibits axis-aligned contour patterns, indicating independent contributions from individual variables. In contrast, CPC-DF22 shows fragmented contour structures caused by interactions among variables, which is characteristic of non-separable function.}
    \label{fig:separability-structures}
\end{figure}

To illustrate this difference, Figure \ref{fig:separability-structures} presents contour visualizations of a separable function (CPC-DF7) and a non-separable function (CPC-DF22).

\subsection{Modality}

Modality of a function refers to the number and spatial distribution of its local and global optima within the search space. Based on modality, functions can be broadly categorized as:

\begin{itemize}
    \item \textbf{Unimodal}: Functions that possess a single global optimum and no other local optima.
    \item \textbf{Multimodal}: Functions that possess multiple local optima, with only one or a subset of them being globally optimal.
\end{itemize}

Compared to unimodal functions, multimodal functions pose substantially greater optimization challenges due to the presence of numerous local optima. This difficulty is particularly pronounced for gradient-based or local search algorithms, which typically converge to the nearest local optimum.

\begin{figure}
    \centering
    \includegraphics[width=0.45\textwidth]{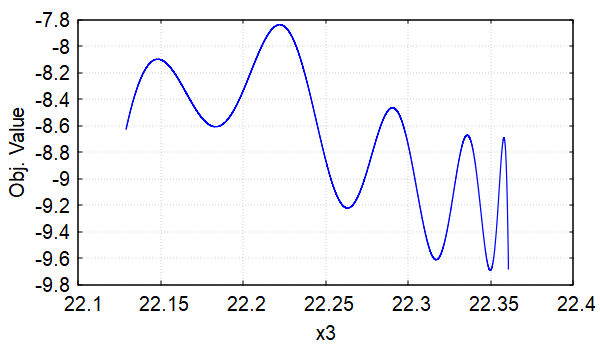}
    \includegraphics[width=0.45\textwidth]{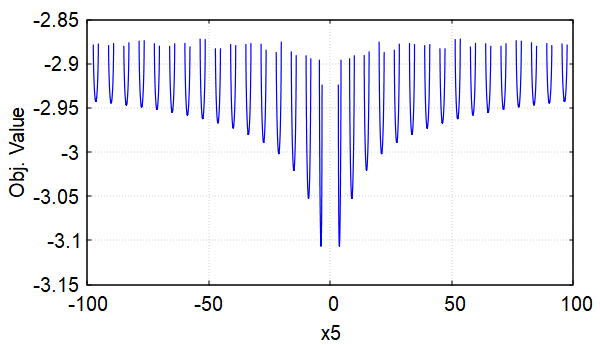}
    \caption{Visualization of multimodal landscapes for CPC-DF6 (left) and CPC-DF12 (right). Multiple distinct peaks and valleys indicate the presence of numerous local optima distributed across disconnected feasible regions.}
    \label{fig:modality-visualization}
\end{figure}

In the CPC benchmark, the functions are highly multimodal and contain numerous local optima distributed across multiple disconnected feasible regions. This fragmented landscape substantially increases the difficulty of global optimization. An effective algorithm must not only avoid premature convergence to local optima within an individual feasible region, but also discover and explore multiple disconnected regions to identify high-quality feasible solution. Figure~\ref{fig:modality-visualization} illustrates the multimodal and fragmented search-space structure commonly observed in the CPC benchmark functions.

\subsection{Multiplicity of Global Optima}

The multiplicity of global optima describes whether an optimization problem admits a single or multiple globally optimal solutions. Based on this property, functions can be classified into the following categories:
\begin{itemize}
    \item \textbf{Single global optimum}: The function admits a unique global optimum.
    \item \textbf{Symmetrical multiple optima}: The function admits multiple global optima due to structural symmetry (e.g., $f(x)=f(-x)$).
    \item \textbf{Asymmetrical multiple optima}: The function admits multiple global optima without structural symmetry.
\end{itemize}

\begin{figure}
    \centering
    \includegraphics[width=0.45\textwidth]{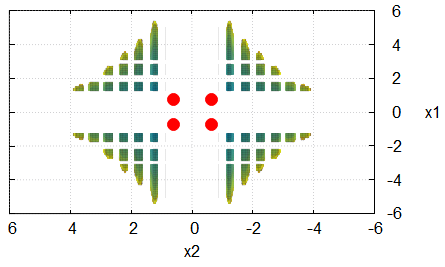}
    \includegraphics[width=0.45\textwidth]{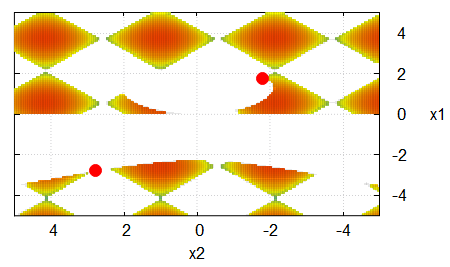}
    \caption{Visualization of global optima structures for CPC-DF5 (left) and CPC-DF12 (right).}
    \label{fig:multiplicity-visualization}
\end{figure}

Examples of symmetrical multiple optima structure and asymmetrical multiple optima structure are presented in Fig.~\ref{fig:multiplicity-visualization}, where the red dots represents the optimal solution. The global optima of CPC-DF5 exhibit clear symmetry, whereas those of CPC-DF12 are distributed asymmetrically in the search space.

To characterize the multiplicity of global optima in the CPC benchmark, we adopt the following counting convention.

If the number of global optima is expressed in the form $2^n$ (e.g., $2^1 = 2$, $2^5 = 32$), the function exhibits symmetrical multiple optima. In such case, identifying any single solution allows all other global optima to be obtained directly through the corresponding symmetry transformations.

If the number is written as a plain integer (e.g., 2, 3) that greater than one, the function contains asymmetrical multiple optima, meaning that the global solutions are not related by symmetry and must be discovered independently. For optimization algorithms, the ability to identify all asymmetrical global optima provides an additional measure of robustness and exploration capability.

\subsection{Solution Sensitivity}

The sensitivity of a solution refers to how small perturbations in the decision variables may lead to significant changes in the objective value or feasibility near the optimal region.

The coexistence of strong nonlinearity and fragmented feasible regions makes the proposed benchmark highly sensitive to such perturbations. Even slight deviations from an optimal solution can produce substantial changes in the objective value or render the solution infeasible. This behavior is particularly evident when the optimum lies near the boundary of the function domain, where the boundary between feasible and infeasible regions can be extremely narrow. Consequently, reliably identifying and remaining near the global optimum becomes considerably more difficult, placing greater demands on the robustness and stability of optimization algorithms.

To illustrate this phenomenon, CPC-DF15 is taken as an example. Table~\ref{tab:cpc-df15-sensitivity} presents the objective values at five points in the vicinity of the best-known optimum of this two-variable function. The boldfaced row corresponds to the best-known optimal solution.

\begin{table}
    \tbl{Sensitivity of objective value near the best-known optimum of CPC-DF15}
    {
        \begin{tabular}{ccc}
            \toprule
            $x_1$ & $x_2$ & Obj. Value \\
            \midrule
            3.17233243403425 & 4.15125018443042 & 9.99800561505972 \\
            \midrule
            \multirow{3}{*}{\textbf{3.17233243403426}} 
            & 4.15125018443041 & 9.99800567819563 \\
            & \textbf{4.15125018443042} & \textbf{9.998005563035544} \\
            & 4.15125018443043 & $\mathrm{NaN}$ \\
            \midrule
            3.17233243403427 & 4.15125018443042 & $\mathrm{NaN}$ \\
            \bottomrule
        \end{tabular}
    }
    \label{tab:cpc-df15-sensitivity}
\end{table}

\begin{figure}
    \centering
    \includegraphics[width=0.45\textwidth]{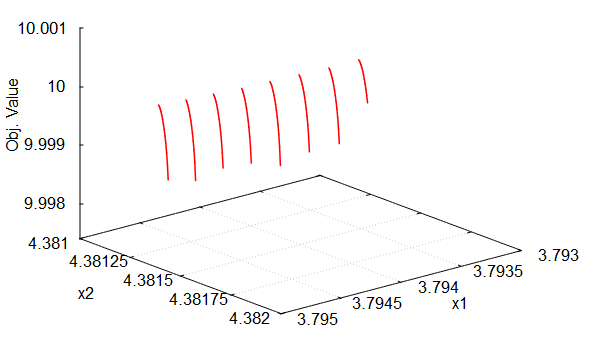}
    \includegraphics[width=0.45\textwidth]{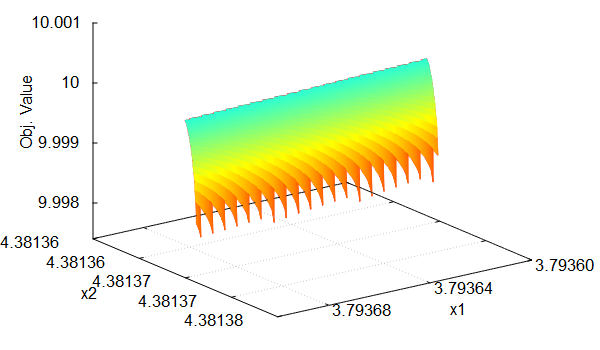}
    \caption{Illustration of the solution sensitivity in CPC-DF15. 
    In the larger range shown on the left, the function appears visually indistinguishable from several straight lines. The magnified view on the right exposes the sharp transitions arising from the vanishingly small gap between feasible and infeasible regions around the optimal solution, which explains the pronounced solution sensitivity.}
    \label{fig:CPC-df15-visualization}
\end{figure}

As illustrated in Fig.~\ref{fig:CPC-df15-visualization}, even a tiny perturbation in $x_1$ or $x_2$ can lead to drastic changes in the objective value and may even render a previously optimal solution infeasible. This behavior highlights the extreme sensitivity of the feasibility boundary in CPC-DF15, and similar patterns are observed across other functions in the benchmark suite. The fact that even infinitesimal perturbations can render a solution infeasible is a defining characteristic of optimization problems with DIDL behavior. It creates highly sensitive regions in the search space where the objective value changes abruptly, thereby imposing substantially higher demands on numerical precision and algorithmic robustness.

This sensitivity is also reflected in our numerical experiments. For several benchmark functions, the objective values reported by some algorithms often deviate noticeably from the values obtained by directly re-evaluating the returned decision variables. In some cases, the reported solutions are even found to be infeasible upon re-evaluation. These discrepancies further demonstrate the extreme sensitivity of optimization problems with DIDL behavior and the stringent precision requirements imposed by hidden feasibility boundaries. They also indicate that accurately identifying feasible and high-quality solutions is substantially more challenging than in conventional continuous optimization settings, as even minor numerical inaccuracies may cause algorithms to misidentify feasibility and move away from the optimal region. A more detailed discussion of these observations is provided in Section~\ref{sec:section-4-numerical-experiment}.

\subsection{Benchmark Functions Overview}

A concise overview of the characteristics of all 25 benchmark functions is presented in Table~\ref{tab:cpc-benchmark-overview}, summarizing these five key properties: dimensionality (Dim.), search range, separability (Sep.), feasible ratio (FR), and the number of global optima (\#Optima).

\begin{table}
    \tbl{Overview of characteristics for CPC Benchmark Functions}
    {
        \begin{tabular}{cccccc}
        \toprule
        Fun ID & Dim. & Search Range & Sep. & FR(\%) & \#Optima \\
        \midrule
        CPC-DF1  & 2      & $[-100,100]^2$    & No  & 48.603052   & 1 \\
        CPC-DF2  & $n=5$  & $[0,10]^n$        & Yes & 0.00114    & 1 \\
        CPC-DF3  & 2      & $[-100,100]^2$    & No  & 0.05724     & 1 \\
        CPC-DF4  & $n=5$  & $[-512,512]^n$    & No  & 0.000484    & 1 \\
        CPC-DF5  & $n=5$  & $[-6,6]^n$        & Yes & 0.000036    & $2^5=32$ \\
        CPC-DF6  & $n=4$  & $[-100,100]^n$    & No  & 0.000544    & 1 \\
        CPC-DF7  & $n=3$  & $[-10,10]^n$      & Yes & 0.007084    & 1 \\
        CPC-DF8  & 2      & $[0,14]^2$        & No  & 15.886572   & 1 \\
        CPC-DF9  & 2      & $[-10,10]^2$      & No  & 37.250012   & 1 \\
        CPC-DF10 & 2      & $[-10,10]^2$      & No  & 35.671264   & 2 \\
        CPC-DF11 & 2      & $[-6,6]^2$        & No  & 0.00114     & $2^2=4$ \\
        CPC-DF12 & $n=5$  & $[-100,100]^n$    & No  & 0.010204    & $2^4=16$ \\
        CPC-DF13 & 2      & $[-100,100]^2$    & No  & 0.000116    & 1 \\
        CPC-DF14 & 2      & $[-10,10]^2$      & No  & 29.208808   & 1 \\
        CPC-DF15 & 2      & $[-5,5]^2$        & No  & 0.000228    & 1 \\
        CPC-DF16 & $n=5$  & $[-512,512]^n$    & No  & 0.00002     & 3 \\
        CPC-DF17 & 2      & $[-100,100]^2$    & No  & 0.554864    & 1 \\
        CPC-DF18 & 4      & $[-5000,5000]^4$  & No  & $\leq 0.000001$ & 1 \\
        CPC-DF19 & 6      & $[-100,100]^6$    & No  & 0.014292    & 1 \\
        CPC-DF20 & $n=5$  & $[-10,10]^n$      & No  & 0.023148    & 1 \\
        CPC-DF21 & $n=5$  & $[-30,30]^n$      & No  & 0.036648    & 1 \\
        CPC-DF22 & 2      & $[-100,100]^2$    & No  & 0.06192     & 1 \\
        CPC-DF23 & $n=5$  & $[-100,100]^n$    & No  & 0.76602     & 1 \\
        CPC-DF24 & 2      & $[-2\pi,2\pi]^2$  & No  & 2.429232    & 2 \\
        CPC-DF25 & $n=3$  & $[-10,10]^n$      & No  & 6.763992    & 1 \\
        \bottomrule
        \end{tabular}
    }
    \label{tab:cpc-benchmark-overview}
\end{table}

\section{Numerical Experiment}\label{sec:section-4-numerical-experiment}

This section presents a comprehensive numerical evaluation of the proposed CPC benchmark using six representative global optimization algorithms. The primary purpose of these experiments is not to establish a definitive ranking among optimization algorithms, but rather to characterize the numerical challenges posed by optimization problems with DIDL behavior and to provide baseline performance references for future studies.

In addition to the CPC benchmark, three representative discontinuity-related benchmark functions from existing benchmark suite are included and evaluated under identical experimental settings. This comparison provides contextual insight into how the proposed benchmark differs from and complements traditional discontinuity-related test problems.

\subsection{Overview of Experimental Design}

All algorithmic settings, computational environments, and evaluation protocols have been described in Section~\ref{subsec:evaluation-protocol}. 

For each benchmark function, all reported solutions are independently re-evaluated using a unified verification procedure. The reported decision variables are substituted back into the objective function to (i) determine whether the solution remains computationally feasible and (ii) compare the re-evaluated objective value with the objective value reported by the algorithm. To reduce the influence of floating-point rounding errors near feasibility boundaries, all verification evaluations are performed using the same double-precision floating-point arithmetic employed in the benchmark implementation. A solution is regarded as feasible only if the re-evaluated objective function returns a valid real-valued result under the verification procedure. The best solution for each algorithm on each benchmark function is selected according to the independently re-evaluated objective values.

To provide a comparative reference, three representative discontinuous functions - F37(Corana Function), F138(Step Function 2), and F170(Xin-She Yang Function 3) - from existing benchmark suites\cite{Jamil2013Survey}, are included in the numerical experiment. These functions respectively represent piecewise structural discontinuity, stepwise quantization discontinuity, and numerically induced pseudo-discontinuity, thereby providing representative reference cases for comparison with the domain-induced discontinuity-like behavior featured in the CPC benchmark.

\subsection{Solution Classification}\label{subsec:solution-classification}

To interpret the numerical behaviors induced by discontinuities and feasibility-sensitive structures, all computed solutions are categorized into four types. This classification framework is essential for distinguishing between valid convergence, local optimality, infeasibility, and numerical inconsistency - phenomena that frequently arise in nonlinear discontinuous optimization. 

\begin{enumerate}
\renewcommand{\labelenumi}{\Roman{enumi}.}
    \item \textbf{Best-Known Solution}: The reported objective value and the corresponding parameters are consistent with the best-known solution.

    \item \textbf{Local Solution}: The obtained solution corresponds to a local solution rather than the best-known solution, either because its objective value is significantly worse or because it lies in a different region of the search space.

    \item \textbf{No Feasible Solution}: The algorithm fails to return any feasible solution within the search domain.

    \item \textbf{Verification-Inconsistent Solution}: The reported objective value is inconsistent with that obtained by independently re-evaluating the reported decision variables using the predefined verification procedure.
    
\end{enumerate}

This classification provides the foundation for analyzing both representative cases and aggregated results across the benchmark. Note that Type IV is not mutually exclusive with Types I–III, as it characterizes the consistency between the reported and independently re-evaluated objective values rather than the optimization quality of the returned solution. In particular, If the reported decision variables become infeasible after re-evaluation, the solution will be classified as Type IV.

\begin{figure}
    \centering
    \includegraphics[width=0.6\linewidth]{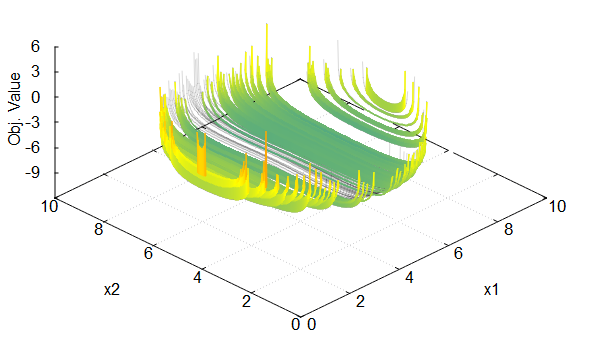}
    \caption{Function graph of CPC-DF3.}
    \label{fig:cpc-demo-df3}
\end{figure}

To better illustrate these four solution types, we consider CPC-DF3 as a representative example. Its analytical form is given in Eq.\eqref{eq:cpc-df3}, and the corresponding function landscape is shown in Fig.~\ref{fig:cpc-demo-df3}. This example highlights the characteristic domain-induced discontinuities that lead to distinct solution behaviors.

\begin{equation}   
    \begin{aligned}
        f(\textbf{x}) := &{- \mathrm{ln} \left( x_{1}-x_{2} \cdot \mathrm{ln} \left( \mathrm{cos} \left( -\frac{2}{3} \cdot x_{1}^{3}-8 \cdot x_{1}^{2} \right) \right) \right) - }\\
        &{\mathrm{ln} \left( 33 \cdot x_{1}-x_{1} \cdot x_{2}+5- \left( \left( x_{1}-4 \right) ^{2}+ \left( x_{2}-5 \right) ^{2}-4 \right) ^{2} \right) }
    \end{aligned}
    \label{eq:cpc-df3}
    \tag{CPC-DF3}
\end{equation}

Where $x_i \in [-100, 100], \quad i = 1, 2$

\begin{table}
    \tbl{Best Solutions and Solution Types on CPC-DF3}
    {\begin{tabular}{ccccc}
        \toprule
        Algorithm & $x_1$ & $x_2$ & Cal.\ Obj.\ Value & Sol.\ Type \\
        \midrule
        Best-Known Solution & \textbf{6.1828121298816} & \textbf{6.490319915658} & \textbf{-10.503674524476093} & -- \\
        UGO  & 6.1828121298816 & 6.49031991565847 & -10.5036630152501 & I, IV \\
        B\&B & 6.5436864435034865 & 4.4179965593342505 & -10.167866753839689 & II \\
        GA   & 6.57657702814721 & 6.16570178994525 & -9.79699243097639 & II \\
        SA   & 6.39262277474369 & 6.34733998281076 & -10.4001096602552 & II, IV \\
        BFGS+MS & \multicolumn{3}{c}{No solution} & III \\
        GPSO & 4.8597692475281562 & 4.5283357354447986 & -8.1619878487073425 & II \\
        \bottomrule
    \end{tabular}}
    \label{tab:cpc-df3-results}
\end{table}

As shown in Table~\ref{tab:cpc-df3-results}, the decision variables returned by UGO coincide with the best-known solution, and the independently re-evaluated objective value also matches the best-known objective value. Therefore, the solution is classified as Type I. However, the objective value reported directly by UGO differs slightly from the independently re-evaluated value, indicating a verification inconsistency. Consequently, the solution is additionally classified as Type IV. This inconsistency affects only the reported objective value and does not alter the optimization quality of the obtained solution.

The solutions obtained by B\&B, GA, SA, and GPSO all yield objective values inferior to the best-known solution and correspond to different locations in the search space. These solutions are therefore classified as Type II. However, the objective value reported by SA is inconsistent with the independently re-evaluated objective value. As a result, the SA result is additionally classified as Type IV.

BFGS+MS fails to return any feasible solution within the prescribed search domain. Consequently, it is classified as Type III.

For Type IV solutions returned by UGO and SA, Table~\ref{tab:cpc-df3-soluiton-type-IV-verification} compares the reported objective values with re-evaluated objective values, revealing the inconsistency.

\begin{table}
    \tbl{Re-evaluated of UGO and SA Solutions for CPC-DF3}
    {
        \begin{tabular}{cccc}
        \toprule
        Algorithm & $x_1$ & $x_2$ & Obj.\ Value \\
        \midrule
        UGO (reported) & 6.1828121298816 & 6.49031991565847 & -10.5036630152501 \\
        UGO (re-evaluated)     & 6.1828121298816 & 6.49031991565847 & -10.503674524476093 \\
        SA (reported) & 6.39262277474369 & 6.34733998281076 & -10.4001096602552 \\
        SA (re-evaluated)     & 6.39262277474369 & 6.34733998281076 & -10.323821540139583 \\
        \bottomrule
        \end{tabular}
    }
    \label{tab:cpc-df3-soluiton-type-IV-verification}
\end{table}

\subsection{Numerical Results}

Building upon the solution classification framework introduced in Section~\ref{subsec:solution-classification}, this section presents the overall numerical performance across all benchmark functions and algorithms.

Table~\ref{tab:cpc-solution-type-distribution} summarizes the distribution of solution types obtained by each algorithm across all 25 benchmark functions, reflecting its overall performance and reliability on optimization problems with DIDL behavior. 

\begin{table}
    \tbl{Distribution of solution types for each algorithm on the 25 CPC benchmark functions}
    {\begin{tabular}{lcccccc}
        \toprule
         & B\&B & UGO & GA & SA & BFGS+MS & GPSO \\
        \midrule
        Type I (\%)   & 12  & 84 & 44 & 0  & 0  & 4  \\
        Type II (\%)  & 64 & 16 & 44 & 68 & 32 & 64 \\
        Type III (\%) & 16 & 0 & 12 & 24 & 68 & 32 \\
        Type IV (\%)  & 20 & 12 & 0  & 12 & 0  & 0  \\
        \bottomrule 
    \end{tabular}}
    \label{tab:cpc-solution-type-distribution}
\end{table}

Table~\ref{tab:cpc-vs-existing-comparison} compares the algorithmic performance on three representative CPC benchmark functions with that on three existing discontinuous benchmark functions. 

\begin{table}
    \tbl{Comparison of algorithmic performance on the CPC benchmark and existing benchmarks}
    {
        \begin{tabular}{lcccccc}
        \toprule
        Function ID & B\&B & UGO & GA & SA & BFGS+MS & GPSO \\
        \midrule
        CPC-DF2  & I  & I & I  & II & III & III \\
        CPC-DF16 & II  & I & III & III & III & III \\
        CPC-DF19 & II, IV & I & II & II & III & II \\
        F37      & I   & I & I  & I  & I   & I  \\
        F138     & I   & I & I  & I  & II  & I  \\
        F170     & I   & I & I  & I  & I   & I  \\
        \bottomrule
        \end{tabular}
    }
    \label{tab:cpc-vs-existing-comparison}
\end{table}

As shown in Table~\ref{tab:cpc-vs-existing-comparison}, the three conventional discontinuity-related benchmark functions are solved successfully by nearly all algorithms, with Type I solutions obtained in almost every case. In contrast, the CPC benchmark frequently produces Type II, Type III, and Type IV outcomes, even for algorithms that perform well on the conventional benchmarks. These observation demonstrate that the CPC benchmark poses substantially greater numerical challenges than the representative discontinuity-related benchmark functions considered in this study.

\section{Discussion}\label{sec:section-5-Discussion}

The numerical experiments demonstrate that the proposed CPC benchmark provides a substantially more challenging test environment than existing discontinuity-related benchmark functions. Rather than serving merely as a platform for comparing optimization algorithms, the experimental results reveal key structural characteristics of optimization problems exhibiting DIDL behavior and demonstrate the benchmark's capacity to differentiate algorithmic behaviors.

Because the feasible regions of many CPC benchmark functions are implicitly defined and highly fragmented, rigorous certification of global optima is generally intractable. Accordingly, all reference solutions reported in this study should be interpreted as best-known solutions rather than mathematically proven global optima. These reference solutions were obtained through extensive computational search using multiple optimization platforms and subsequently verified through an independent evaluation procedure.

It should also be noted that the optimization algorithms were implemented in different software environments, as introduced in Section ~\ref{subsec:overview-of-selected-algorithms}, each employing different numerical mechanisms, infeasibility-handling strategies, and stopping criteria. Consequently, the observed performance differences reflect not only algorithmic characteristics but also implementation-dependent numerical behavior. Therefore, the reported results should be interpreted primarily as baseline performance references for the proposed benchmark rather than as a definitive ranking of optimization algorithms.

\subsection{Comparative Analysis of Discontinuity-Related Benchmarks} \label{subsec:comparative-analysis-of-discontinuous-related-benchmarks}

The numerical results in Section~\ref{sec:section-4-numerical-experiment} reveal a clear contrast between the CPC benchmark and existing benchmarks. For the representative discontinuity-related benchmark functions F37, F138, and F170, all algorithms except BFGS+MS consistently achieve Type I solutions, indicating that these benchmarks pose limited difficulty for modern global optimization methods. Their structural simplicity - whether arising from piecewise definition, stepwise quantization, or pseudo-discontinuity - results in landscapes that are easily navigated and provide limited discrimination between optimization algorithms.

In contrast, the CPC benchmark presents markedly greater difficulty. For functions such as CPC-DF2, CPC-DF16, and CPC-DF19, algorithms frequently fail to reach the best-known solution, with many runs instead yielding Type II, Type III, or Type IV outcomes. These patterns reflect the irregular feasible regions and implicit structural constraints inherent to the CPC functions, which give rise to sharp boundary effects and highly nonlinear local landscapes.

The divergence in algorithmic behavior across CPC functions demonstrates that the benchmark provides a richer and more discriminative test environment. Algorithms that perform similarly on traditional discontinuous benchmarks exhibit clear performance separation on the CPC benchmark, revealing differences in robustness, feasibility-handling mechanisms, and sensitivity to numerical precision. Overall, these findings indicate that existing discontinuous benchmarks are insufficient for assessing algorithmic performance under realistic structural discontinuities, whereas the CPC benchmark offers a more challenging and informative testbed.

\subsection{Algorithmic Behavior under Domain-Induced Discontinuity-Like Structures}

\subsubsection{Limitations of Gradient-Based Algorithms}

Domain-induced discontinuity-like structure and fragmented feasible regions fundamentally undermine the assumptions required by gradient-based optimization methods. In such landscapes, gradient information becomes unreliable or undefined across large portions of the search space, and disconnected feasible regions prevent gradient directions from guiding the search toward global optima. Even with reasonably chosen initial points, these structural properties limit the ability of gradient-dependent methods to escape infeasible gaps or navigate toward distant feasible regions.

\begin{figure}
    \centering
    \includegraphics[width=0.45\textwidth]{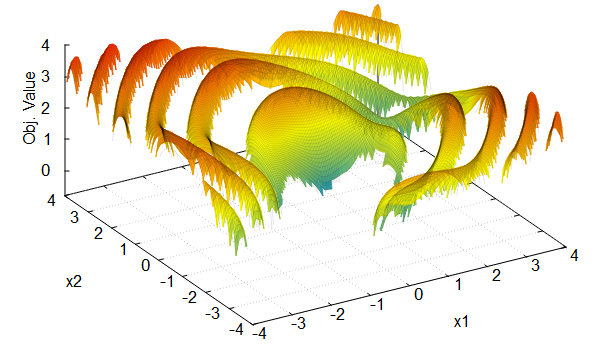}
    \includegraphics[width=0.45\textwidth]{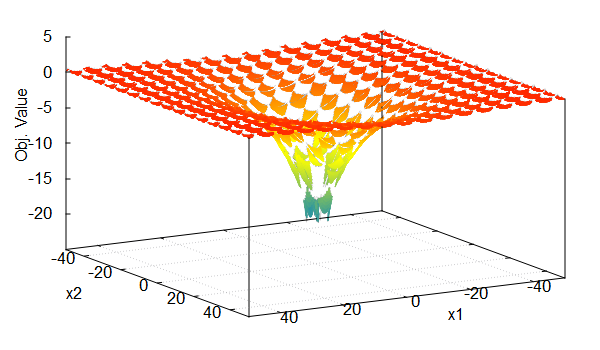}
    \caption{Illustration of structural discontinuity-like behavior and fragmented feasible regions in the CPC Benchmark (CPC-DF1 and CPC-DF23).}
    \label{fig:demo for discontinuity}
\end{figure}

As shown in Table.~\ref{tab:cpc-solution-type-distribution}, BFGS+MS fails to produce any Type I solution, and GPSO yields Type II or Type III outcomes in 88\% of cases. Figure~\ref{fig:demo for discontinuity} illustrates representative examples of CPC benchmark functions with disconnected feasible regions where gradient-based methods tend to converge prematurely or fail to identify feasible regions due to the lack of reliable directional information.

These observations highlight that the difficulty is intrinsic to the landscape rather than to algorithmic design. Gradient-based methods are inherently limited when applied to problems with DIDL behavior, as they cannot rely on smoothness or continuity to guide the search and must instead operate under incomplete or misleading derivative information.

\subsubsection{Impact of Feasibility Scarcity and a Practical Mitigation Strategy}

A defining characteristic of the CPC benchmark is the extremely low feasible ratio (FR) exhibited by many functions. Table~\ref{tab:cpc-fr-group-distribution} shows a pronounced shift in solution types when FR falls below 1\%, with Type III solutions increasing from 4.8\% to 31.5\%. When feasible regions occupy only a tiny fraction of the search space, the primary challenge shifts from optimization itself to the discovery of feasible regions, and algorithms may fail to identify any feasible solution throughout the entire run. 

To mitigate this difficulty, a practical strategy is to incorporate a Monte-Carlo-based preprocessing phase that samples a large number of candidate points before optimization begins. This increases the likelihood of locating at least one feasible point for initialization and therefore provides a practical means of improving optimization performance in low-feasibility settings. While this does not eliminate the underlying structural difficulty, it provides a simple and effective mechanism for enhancing algorithmic stability on CPC functions with extremely scarce feasible regions.

\begin{table}
    \tbl{Solution type distribution for CPC benchmark functions grouped by feasible ratio (FR) above and below 1\%}
    {
        \begin{tabular}{lcccc}
        \toprule
        Function type & Type I (\%) & Type II (\%) & Type III (\%) & Type IV (\%) \\
        \midrule
        FR $> 1\%$  & 26.2 & 64.3 & 4.8  & 9.6 \\
        FR $< 1\%$  & 26.9 & 38.0 & 34.3 & 4.6  \\
        \bottomrule
        \end{tabular}
    }
    \tabnote{Among the 25 CPC benchmark functions, seven (CPC-DF1, CPC-DF8, CPC-DF9, CPC-DF10, CPC-DF14, CPC-DF24 and CPC-DF25) have a feasible ratio (FR) greater than 1\%, while the remaining functions have FR below 1\%.}
    \label{tab:cpc-fr-group-distribution}
\end{table}

\subsubsection{Numerical Sensitivity near Feasibility Boundaries}

Domain-induced discontinuity-like behavior inherently creates extremely sharp and irregular feasibility boundaries. Unlike conventional benchmark functions, feasibility in CPC functions is embedded implicitly within nonlinear objective expressions rather than specified by explicit constraints. As a result, the best-known solutions of several CPC functions lie very close to the feasibility boundary, where minute perturbations of the decision variables may abruptly invalidate the objective function.

CPC-DF3 provides a representative example of this phenomenon. As shown in Table~\ref{tab:cpc-further-df3-validation}, two neighboring solutions obtained by perturbing the last decimal digit of the best-known solution are examined.

\begin{table}
    \tbl{Objective values of two neighboring solutions around the best-known solution of CPC-DF3}
    {
        \begin{tabular}{cccc}
        \toprule
        Solution & $x_1$ & $x_2$ & Obj. Value\\
        \midrule
        A (slight below) & 6.1828121298815	& 6.49031991565847 & -10.503674524476093\\
        B (Best-known) & \textbf{6.1828121298816} & \textbf{6.49031991565847} & \textbf{-10.192474818519308}\\
        C (slightly above)& 6.1828121298817	& 6.49031991565847 & $\mathrm{NaN}$\\
        \bottomrule
        \end{tabular}
    }
    \label{tab:cpc-further-df3-validation}
\end{table}

The objective function of CPC-DF3 is defined as:

\begin{equation}   
    \begin{aligned}
        f(\textbf{x}) := &{- \mathrm{ln} \left( x_{1}-x_{2} \cdot \mathrm{ln} \left( \mathrm{cos} \left( -\frac{2}{3} \cdot x_{1}^{3}-8 \cdot x_{1}^{2} \right) \right) \right) - }\\
        &{\mathrm{ln} \left( 33 \cdot x_{1}-x_{1} \cdot x_{2}+5- \left( \left( x_{1}-4 \right) ^{2}+ \left( x_{2}-5 \right) ^{2}-4 \right) ^{2} \right) }
    \end{aligned}
    \label{eq:cpc-df3}
    \tag{CPC-DF3}
\end{equation}

Therefore, the objective function of CPC-DF3 induces three implicit feasibility conditions expressed as follows.

\begin{equation}   
    \left\{
        \begin{aligned}
            &base_{3\text{-}1} := \mathrm{cos} \left( -\frac{2}{3} \cdot x_{1}^{3}-8 \cdot x_{1}^{2} \right) > 0\\
            &base_{3\text{-}2} := x_{1}-x_{2} \cdot \mathrm{ln} \left( base_{3\text{-}1} \right) > 0\\
            &base_{3\text{-}3} := 33 \cdot x_{1}-x_{1} \cdot x_{2}+5- \left( \left( x_{1}-4 \right) ^{2}+ \left( x_{2}-5 \right) ^{2}-4 \right) ^{2} > 0
        \end{aligned}
    \right.
    \label{eq:cpc-df3-feasibility condition}
    \tag{Feasibility Condition for CPC-DF3}
\end{equation}

\begin{table}
    \tbl{Feasibility Verification for Solution B and C}
    {
        \begin{tabular}{cccc}
        \toprule
        Solution & Variable & Value & Condition \\
        \midrule
        B & $base_{3\text{-}1}$ & 1.4763849459792957e-15 & $> 0$ \\
        B & $base_{3\text{-}2}$ & 227.82191453323668 & $> 0 $ \\
        B &$base_{3\text{-}3}$ & 159.989834258271 & $> 0$ \\
        C & $base_{3\text{-}1}$ & -1.7449453205322083e-11 & $< 0$ \\
        C & $base_{3\text{-}2}$ & Cannot compute & - \\
        C &$base_{3\text{-}3}$ & 159.98983425827097 & $> 0$ \\
        \bottomrule
        \end{tabular}
    }
    \label{tab:cpc-df3-feasibility-check}
\end{table}

Table~\ref{tab:cpc-df3-feasibility-check} shows that Solution B satisfies all three implicit conditions, whereas Solution C violates the first condition because the argument of logarithm must be positive, rendering the logarithmic term undefined. Although Solution C differs from Solution B only in the final decimal digit of $x_1$, the perturbation changes the sign of $base_{3\text{-}1}$ from positive to negative, violating the implicit feasibility constraint, rendering the logarithmic term undefined and causing the objective function to return $\mathrm{NaN}$.

The example above illustrates the extreme sensitivity of the feasible region in CPC-DF3 to minute perturbations of the decision variables. A change in only the last decimal digit of the best-known solution is sufficient to alter the feasibility status, causing the objective function to become undefined. This behavior is not a consequence of implementation-dependent numerical instability, but instead arises from the intrinsic domain restrictions embedded in the analytical form of the objective function. Similar situations are observed in several CPC benchmark functions, where the best-known solutions lie extremely close to implicit feasibility boundaries. Consequently, optimization algorithms must simultaneously search for high-quality solutions while preserving feasibility throughout the optimization process. This characteristic distinguishes the CPC benchmark from conventional benchmark suites, where objective functions are typically well-defined throughout the prescribed search domain and feasibility rarely becomes a dominant source of optimization difficulty.

\subsubsection{Domain-Induced Discontinuity-Like Behavior as Hidden-Constrained Optimization}

Domain-induced discontinuity-like behavior can be interpreted as a form of implicit feasibility-constrained optimization. Although no explicit constraint functions are specified, meaningful objective evaluation is only possible within the natural domain of the objective function. Consequently, the optimization process is effectively constrained by hidden feasibility conditions embedded in the analytical expression itself. In other words, the optimization algorithm is required not only to search for better objective values but also to remain within an implicitly defined feasible region throughout the search.

CPC-DF3 provides a representative example of this implicit-feasibility structure. Meaningful objective evaluation requires satisfying the conditions $base_{3-1} > 0$, $base_{3-2} > 0$, and $base_{3-3} > 0$, although none of these conditions are explicitly supplied to the optimization algorithm. Unlike conventional constrained optimization problems, where constraint functions are explicitly available and feasibility can be evaluated independently, the CPC benchmark provides only the objective function. The implicit feasibility conditions arise naturally from mathematical operators such as logarithms, fractional powers, and exponential expressions, rather than from separately defined constraint equations. Consequently, optimization algorithms must infer the admissible region solely through objective evaluations and their responses to undefined function values.

From an algorithmic perspective, this interpretation fundamentally changes the optimization task. The primary challenge is no longer limited to objective minimization; algorithms must simultaneously discover the hidden feasible region while optimizing the objective. This requires reliable mechanisms for detecting infeasible evaluations, maintaining numerical stability near feasibility boundaries, and recovering from invalid function evaluations without interrupting the search.

This interpretation also clarifies why the CPC benchmark differs fundamentally from conventional discontinuous benchmark functions. Rather than evaluating an algorithm solely on its ability to optimize a discontinuous objective landscape, the CPC benchmark simultaneously evaluates its capability to discover and maintain feasibility in the presence of hidden domain restrictions. Consequently, the benchmark places greater emphasis on feasibility exploration, numerical robustness, and boundary-aware search, all of which are common characteristics of practical nonlinear optimization problems.

\subsubsection{Challenges in Explicitly Reformulating Implicit Feasibility Conditions}

A seemingly straightforward way to address domain-induced discontinuities is to extract the implicit feasibility conditions from the objective function and rewrite them as explicit constraints. In some CPC functions, this approach is indeed applicable and can simplify the optimization process by preventing undefined evaluations. For example, in CPC-DF13 from Section~\ref{subsec:evaluation-protocol}, the only two implicit conditions $base_{13\text{-}1} \ge 0$ and $base_{13\text{-}2} \ge 0$ can be stated explicitly, transforming the hidden feasibility requirements into verifiable constraints and reducing the likelihood of encountering undefined objective values. This can increase the chance that an algorithm successfully identifies feasible points and reaches to a better solution. Nevertheless, the ability to solve CPC functions without relying on such explicit extraction reflects a stronger level of algorithmic robustness, as the algorithm must infer feasibility solely from objective evaluations.

However, explicitly extracting these feasibility conditions is not universally applicable. In practice, manually deriving all implicit feasibility conditions can itself become cumbersome and error-prone, particularly when the admissible domain cannot be expressed as a single fixed condition, but instead changes according to the values of the decision variables. For real-valued objective evaluations in standard numerical optimization environments, expressions of the form $x_1^{x_2}$ provide a representative example, whose admissible domain can be characterized as follows:

A representative case in CPC benchmark is CPC-DF14, whose definition is given below:

\begin{equation*}   
    \begin{aligned}
        f(\textbf{x}) := &{ \left( 1.5-x_{1}+x_{1} \cdot x_{2} \right) ^{0.2}+ }\\
        &{\left( 2.5-x_{1}+x_{1} \cdot x_{2}^{2} \right) ^{x_{1}-0.5}+ \left( 2.625-x_{1}+x_{1} \cdot x_{2}^{3} \right) ^{0.2}}
    \end{aligned}
    \label{eq:cpc-df14}
    \tag{CPC-DF14}
\end{equation*}

Where $x_i \in [-10, 10], \quad i = 1, 2$.

To illustrate this difficulty, we attempt to explicitly derive all implicit feasibility conditions for CPC-DF14. The three terms that introduce implicit domain restrictions are listed below.

\begin{equation}   
    \left\{
        \begin{aligned}
            &base_{14\text{-}1} := 1.5-x_{1}+x_{1} \cdot x_{2}&, \; base_{14\text{-}1} \geq 0\\
            &base_{14\text{-}2} := 2.5-x_{1}+x_{1} \cdot x_{2}^{2}&, \; (base_{14\text{-}2})^{x_1-0.5} \; is \; well\text{-}deinfed\\
            &base_{14\text{-}3} := 2.625-x_{1}+x_{1} \cdot x_{2}^{3}&, \; base_{14\text{-}3} \geq 0
        \end{aligned}
    \right.
    \label{eq:cpc-df3-feasibility condition}
    \tag{Feasibility Condition for CPC-DF14}
\end{equation}

The feasibility conditions associated with $base_{14\text{-}1}$ and $base_{14\text{-}3}$ are straightforward, since both terms simply require non-negative bases. By contrast, the middle term $base_{14\text{-}2}$ is considerably more complicated because its feasibility condition itself is case-dependent rather than fixed.

\begin{equation*}  
    (base_{14\text{-}2})^{x_1-0.5} \quad well\text{-}defined \iff
        \left\{
            \begin{aligned}
                &if \; x_1-0.5 \in \mathbb{Z}_{\geq 0}, then \; base_{14\text{-}2} \in \mathbb{R}\\
                &if \; x_1-0.5 \in \mathbb{Z}_{< 0}, then \; base_{14\text{-}2} \in \mathbb{R} \setminus \{0\}\\
                &if \; x_1-0.5 \in (\mathbb{R} \setminus \mathbb{Z})_{\ge 0}, then \; base_{14\text{-}2} \geq 0\\
                &if \; x_1-0.5 \in (\mathbb{R} \setminus \mathbb{Z})_{< 0}, then \; base_{14\text{-}2} > 0            \end{aligned}
        \right.
    \tag{Feasibility Condition for $base_{14\text{-}2}$}
\end{equation*}

Since the exponent $x_1 - 0.5$ varies with the decision variable, different feasibility conditions apply in different regions of the search space. Consequently, the admissible domain of the middle term cannot be represented by a single explicit constraint but instead requires a case-by-case decomposition. This example illustrates that even for a two-dimensional objective function, manually reformulating all implicit feasibility conditions may already become cumbersome.

A common adopted strategy for handling the expressions of $x^y$ is to apply the transformation $x^y = exp \left( y\ln(x) \right)$, an approach employed by several state-of-the-art global optimization solvers\cite{BARONmanual2026, GAMSuserguide}. Under the computationally admissible domain adopted in this study, however, this reformulation is not equivalent to the original power expression because it implicitly assumes a strictly positive base. As a result, valid evaluations permitted by the original expression when the exponent is an integer are excluded, thereby altering the computationally admissible domain and potentially modifying the optimization problem itself.

In summary, explicitly reformulating implicit feasibility conditions can, in some cases, reduce the occurrence of invalid objective evaluations and improve the probability of locating feasible solutions. However, this approach is not universally applicable. While some objective functions admit fixed feasibility conditions that can be extracted directly, others involve case-dependent admissible domains, where the applicable feasibility condition vary with the values of the decision variables. Furthermore, commonly used reformulations of nonlinear expressions may alter the computationally admissible domain and therefore fail to preserve the original optimization problem. Consequently, the ability to solve CPC function without relying on manually extracted feasibility conditions represents a stronger level of algorithmic robustness, particularly in black-box optimization where the underlying feasibility structure is unavailable.

\section{Conclusion and Future Work}\label{sec:section-6-conclusion}

This paper presented the CPC benchmark, a collection of 25 nonlinear optimization problems exhibiting domain-induced discontinuity-like (DIDL) behavior arising from intrinsic domain restrictions embedded within objective functions. Unlike conventional discontinuous benchmarks, the CPC benchmark characterizes optimization problems whose feasible regions are determined implicitly by the computationally admissible domain rather than by explicit constraint functions.

The study identifies several structural characteristics that distinguish DIDL problems from conventional benchmark functions, including hidden feasibility structures, extremely low feasible ratios, fragmented and disconnected feasible regions, and pronounced numerical sensitivity near feasibility boundaries. Together, these characteristics transform optimization into a task that requires not only objective optimization but also implicit feasibility discovery.

The numerical experiments indicate that the CPC benchmark presents substantially greater challenges than representative conventional discontinuity-related benchmark functions. Many algorithms fail to reach the best-known solutions, while methods exhibiting similar performance on existing benchmarks become clearly differentiated on the CPC benchmark. These observations suggest that the CPC benchmark provides a more informative evaluation environment for assessing robustness, feasibility handling, and numerical stability in optimization problems involving hidden feasibility structures.

For future work, we plan to further expend the CPC benchmark framework to encompass a broader spectrum of challenging optimization problems, including systems of nonlinear equations, curve fitting problems, and ill-conditioned optimization problems. All current and future benchmark instances, together with benchmark documentation and reference implementations will be made publicly available through the CPC Benchmark GitHub repository. The repository is intended to serve as a comprehensive and extensible platform for the systematic evaluation and development of advanced global optimization algorithms. Details of the repository are provided in the Software Availability section.

We hope that the CPC benchmark will encourage the development of optimization algorithms capable of handling hidden feasibility structures that frequently arise in practical scientific and engineering applications but remain underrepresented in existing benchmark collections.

\section*{Acknowledgment(s)}
I would like to thank Professor Makoto Yamashita for his valuable guidance and insightful suggestions throughout the development of this benchmark study.

\section*{Disclosure statement}
The authors report no conflict of interest.

\section*{Funding}
This research received no specific grant from any funding agency.

\section*{Software Availability}
The scripts used for the numerical experiments, the benchmark implementation, are available at \href{https://github.com/CPC-Opt/CPC-Benchmark}{https://github.com/CPC-Opt/CPC-Benchmark}, while MATLAB R2024b, PlatEMO, LINGO 20.0, and 1stOpt 11.0 were used under their respective license terms.

\section*{Data Availability}
The definitions of the 25 benchmark functions, their search domains, the reference solutions used in this paper are available at \href{https://github.com/CPC-Opt/CPC-Benchmark}{https://github.com/CPC-Opt/CPC-Benchmark}
.

\section*{Reproducibility Statement}
All experiments were conducted on a Windows 11 workstation with an Intel Core i7-12700H CPU and 32 GB RAM using double-precision arithmetic, and the software versions, benchmark definitions are provided in the manuscript and repository.

\end{document}